\documentclass[12pt]{article}

\usepackage{amsmath,amssymb,amscd,graphicx,epsfig}

\everymath{\displaystyle}

\newtheorem{teo}{Theorem}[section]
\newtheorem{lem}[teo]{Lemma}
\newtheorem{cor}[teo]{Corollary}

\newtheorem{prop}[teo]{Proposition}
\newtheorem{lem-defi}[teo]{Lemma-Definition}

\newtheorem{remark}[teo]{Remark}

\newcommand{\mt}{\mathbb{T}}
\newcommand{\mr}{\mathbb{R}}
\newcommand{\me}{\mathbb{E}}
\newcommand{\mg}{\mathbb{G}}

\newcommand{\mz}{\mathbb{Z}}

\newcommand{\Css}{$C^{\ast}$-algebras }           

\def\B{{\cal B}}

\def\F{{\cal F}}

\def\C{{{\cal C}}}

\def\L{{{\cal L}}}

\def\P{{{\cal P}}}
\def\S{{{\cal S}}}

\def\V{{{\cal V}}}

\def\M{{{\cal M}}}

\def\F{{\cal F}}

\def\C{{{\cal C}}}

\def\L{{{\cal L}}}

\def\P{{{\cal P}}}
\def\S{{{\cal S}}}

\def\V{{{\cal V}}}

\title{On the Dynamics of $\mg$-Solenoids\\
Applications to Delone Sets
}

\author{\normalsize  Riccardo Benedetti* \,and 
Jean-Marc
Gambaudo**}

\begin{document}

\maketitle

\begin{abstract} 
 A  $\mg$-solenoid is a laminated space whose leaves are copies of a single Lie group $\mg$,  and whose transversals are totally disconnected sets. It inherits a  $\mg$-action and can be considered as dynamical system. Free $\mz^d$-actions on the Cantor set as well as a large class of tiling spaces possess  such a structure of $\mg$-solenoid. We show that a  $\mg$-solenoid can be seen as a projective limit of branched manifolds modeled on $\mg$ . This allows us to give a topological description of the transverse invariant measures associated with a $\mg$-solenoid in terms of  a positive cone in the projective limit of the $dim (\mg)$-homology groups of these branched manifolds. In particular we exhibit a simple criterion implying unique ergodicity.  A particular attention is paid to the case when the Lie group $G$ is the group of affine orientation preserving isometries of the Euclidean space or its subgroup of translations.

 \end{abstract}
\bigskip
\bigskip
\bigskip

\bigskip
\bigskip
\bigskip
\noindent \small{(*) Dipartimento di
Matematica, via F. Buonarroti 2, 56127 Pisa, Italy.}

\noindent \small{(**) Laboratoire de Topologie, U.M.R. 5584 du CNRS, 
Universit\'e de Bourgogne,
B.P. 47870- 21078 Dijon Cedex France. }

\eject

\section{ $\mg$-solenoids}

  Let $\mg$ be a connected Lie group and $d_\mg$ a right invariant metric on $\mg$. Consider  a compact metric space $M$ and assume there exist a cover
of $M$ by open sets $U_i$ called {\it boxes} and homeomorphisms called {\it charts}
$h_i:U_i\to V_i\times T_i$ where the $T_i$'s are totally disconnected metric sets and the $V_i$'s are open subsets in   $\mg$.  These open sets and
homeomorphisms define an atlas of a {\it $\mg$-solenoid}
structure on $M$ if the {\it
transition maps} $h_{i,j}\, =\, h_j\circ h_i^{-1}$ read on their
domains of definitions: $$h_{i, j}(v, t)\,=\, (v. g_{i,j},\,
\tau_{i, j}(t)),$$ where $\tau_{i, j}$ is  a continuous map, and $g_{i,j}$ is an  element  in $\mg$.
 Two atlas are {\it
equivalent} if their union is again an atlas. The fact that the transition maps are required to be very rigid implies the following  properties:
\begin{itemize}
\item [(1)] {\it Laminated structure}. We call {\it slice} of a solenoid,  a subset of the form
$h_i^{-1}(V_i\times \{t\})$.  The transition maps map slices onto slices.  The {\it leaves} of $M$
are the smallest connected subsets that contain all the slices
they intersect. 
 Each leaf is a  differentiable manifold with dimension  $dim(\mg)$. 

\item [(2)] { \it   Riemannian metric}. The right invariant metric $d_\mg$ on $\mg$ induces a Riemannian metric $d_L$ on each leaf $L$ in $M$.
\item [(3)] {\it Tranverse germ}. We call {\it vertical} of a solenoid, a subset of the form $h_i^{-1}(\{x\}\times T_i)$. The 
transition maps map verticals onto verticals.  This allows us to define above each point in the $\mg$-solenoid 
a local vertical (independently on the charts).
\end{itemize}
A {\it $\mg$-solenoid}
is the data of a metric compact space $M$ together
with an equivalence class of atlas which satisfies the further hypothesis:

\noindent Each leaf $L$ in $M$ equipped with is metric $d_L$ is isometric to  $\mg$ equipped with its metric $d_\mg$.

$\mg$-solenoids satisfy the following extra properties:
\begin{itemize}
\item [(4)]  { \it $\mg$-action}.  The free right action on $\mg$  induces a free  right $\mg$-action on $M$  whose orbits are the leaves of $M$. This action is   free on the leaves $L$  isometric to $\mg$.

\item [(5)] {\it Locally constant return times}. Since the $\mg$-action preserves the verticals, it follows that if  a leaf $L$ ( the $\mg$-orbit of some point) intersects two verticals $V$ and $V'$ at $v$ and $v.g$ where $g\in\mg$ then,  for any $\hat v$ in $V$ close enough to $v$, $\hat v .g $ is a point in $V'$ close to $v'$.\end{itemize} 

A   $\mg$-solenoid $M$ is {\it  minimal} if each leaf of $M$ is dense in $M$. It is {\it  expansive} if the $\mg$-action is expansive {\em  i.e.} there exists a positive constant $\epsilon (M)$ such that if two points on the same local vertical have their $\mg$-orbits which remain $\epsilon(M)$-close one of the other then these two points coincide. The constant $\epsilon (M)$ is called {\it the expansivity  diameter} of the $\mg$-solenoid.

\bigskip
\noindent { \sl Standard examples}

\begin{enumerate}

\item The most naive example of  $\mg$-solenoid is given by a compact, connected Lie group. In this case all the verticals are reduced to a single point.
\item More interesting is the case of the dynamical system given by the iteration of an homeomorphism  $\phi$ on the Cantor set $\Sigma$ with no periodic point. It is plain to see that the suspension $(\mt^1\times \Sigma)/ \phi$ inherits a structure of  $\mr$-solenoid.  On one hand, the standard odometer is an  example of non expansive  $\mr$-solenoid, on the other hand, the shift operation on a minimal non periodic subset of $\{0, 1\}^\mz$ is an expansive  $\mr$-solenoid.
\item More generally any free $\mz^d$-action on the Cantor set inherits a structure of  $\mr^d$-solenoid. 
\end{enumerate}

\noindent This paper is organized as follows:

\bigskip \noindent $\bullet$ Section 2 is devoted to provide a much broader set of examples of  $\mg$-solenoids. Let $\Gamma$ be the maximal compact subgroup in $\mg$. The quotient space    $\mg/\Gamma$ is a simply connected manifold which inherits  a Riemannian metric such that $\mg$ acts transitively by isometries on $\mg/\Gamma$.  On this Riemannian manifold one can define Delone sets of  finite $\mg$-type. The $\mg$-action on such  a set $X$ can be extended to  a $\mg$-action on the continuous hull $\Omega_\mg(X)$. We prove that, when the Delone set is  strongly aperiodic , this dynamical system has a natural structure of expansive $\mg$-solenoid (Theorem \ref{sol0}). Conversely, we prove that  any expansive and minimal $\mg$-solenoid can be seen as a continuous hull of a strongly aperiodic Delone set of finite $\mg$-type (Theorem \ref{sol}). 

\noindent $\bullet$ In section 3, we develop a combinatorial approach of the structure of  $\mg$-solenoid (Theorem \ref{tour}). This approach can be interpreted as an extension to our general situation of the construction of infinite sequences of towers and Brateli diagrams that turned to be so powerful in the study of $\mz$-action on the Cantor set.

\noindent $\bullet$ This allows us to introduce in Section 4,  a special class of branched manifolds modeled on $\mg$, that we call $\mg$-branched manifolds and to show   that a  $\mg$-solenoid is  a projective limit of $\mg$-branched manifolds whose faces become arbitrarily large when going backward in the projective limit (Theorem \ref{branched}). This construction gives us better and better approximants of  $M$ with bigger and bigger regions where the $\mg$-action is defined.  The known analogue of this description in the case of the $\mz$-action by a homeomorphism $\phi$ on the Cantor set $\Sigma$, is that the dynamical system $(\Sigma, \phi)$ can be seen as a projective limit of oriented  graphs with  all vertices but one  with valence 2, and such that the length of the loops around the singular vertex goes to infinity (see for instance \cite{Vershik},  \cite{Putnam},\cite{GM}). In the context of smooth dynamics, projective limit of branched manifolds where introduced by R. Williams in the $60$'s  (\cite {W1}, \cite {W2}). Our approach makes a bridge between these last two  constructions.

\noindent $\bullet$ This projective limit construction allows us to give in Section 5, a topological description of  the set of tranverse invariant measures of a  $\mg$-solenoid as the projective limit of cones in the $dim\, \mg$-homology groups of the $\mg$-branched manifolds (Theorem \ref{measure}). We show that if the number of faces of the branched manifolds is uniformly bounded then  this number is also a  bound for the number of ergodic transverse invariant measures (up to a scaling factor). Furthermore, we prove that if the linear maps projecting the homology groups one into the other in the projective limit, are bounded then there exists a unique transverse invariant measure (up to a scaling factor) (Corollary \ref{cor}).

\noindent $\bullet$ Finally, Section 6 is devoted to a more specific analysis of the case when the Lie group is the group $\me^d$ of isometries in $\mr^d$ or the group of translations $\mr^d$.

\bigskip
\begin{remark} {\rm Our definition of a  $\mg$-solenoid implies that the $\mg$-action its wears is free. We could have worked out a same analysis for the case when some leaves of the $\mg$-solenoid are  quotients of $\mg$ by a finite group. This would have allowed us deal also to with  Delone sets which are invariant under finite order isometries (see Remark \ref{re}).  We made this choice only for reasons of  clearness, in order to work only with Riemannian manifolds and not orbifolds. Carrying on all along our discussion the this distinction "manifold-orbifold" would have been too cumbersome. However the whole construction of projective  limits works in the same way as well as the topological interpretation of the transverse invariant measure.}
\end{remark}

\bigskip

\begin{remark} {\rm  For $\mr^d$-action, J. E. Anderson and I. Putnam \cite{AP} proved that the hull $\Omega_{\mr^d}(Y)$ of a  Delone set with finite $\mr^d$-type obtained by a substitution rule is a projective limit of branched manifolds  whose faces become arbitrarily large when going backward in the projective limit. This result has been extended by J. Bellissard,  R. Benedetti and  J.-M. Gambaudo \cite{BBG} to the general case of  a  Delone set with finite $\mr^d$-type. F. G\" ahler \cite{Gahler} announced recently a very simple  proof of the fact that the  hull is a projective limit of branched manifolds in the same context. However in this proof, the sizes of the faces of the branched manifolds remain constant as we go backward in the projective limit and thus, the number of faces goes to infinity. More recently, L. Sadun has extended the construction of F. G\" ahler to the hull $\Omega_{\mg}(Y)$ of a  Delone set with finite $\mg$-type $Y$ in a homogeneous Riemannian manifold \cite{Sad}.  
Finally the fact that the hull $\Omega_{\mg}(Y)$ of a  Delone set with finite $\mg$-type $Y$ in a homogeneous Riemannian manifold has a laminated structure was  observed by E. Ghys \cite{Ghys}.}
\end{remark}

\section {Delone sets: the dynamical point of view}

 Let $\mg$ be a connected Lie group, $\Gamma$ its maximal compact subgroup. The quotient space    $\mg/\Gamma$ (right classes)  is a simply connected manifold which inherits  a Riemannian metric such that $\mg$ has a transitive right action  by isometries on $\mg/\Gamma$.  

\noindent Consider two positive numbers $r$ and $R$, an {\it $(r,R)$-Delone set} is a discrete subset $X$ of $\mg/\Gamma$   which satisfies the following two properties:
\begin{itemize}
\item [(i)] {\it Uniformly Discrete.} Each open ball with radius $r$ in $\mg/\Gamma$ contains at most one point in $X$;
\item[(ii)] {\it Relatively Dense.} Each open ball with radius $R$ contains at least one point in $X$.
\end{itemize}
When  the constants $r$ and $R$ are not explicitely used, we will say in short {\it Delone set} for an $(r, R)$-Delone set. We refer to  \cite{LP} for a more detailed approach of the theory of Delone sets. 

 A {\it patch} of a Delone set is a finite subset of $X$. Two patches $P_1$ and $P_2$  of  a Delone set $X$ are of the same $\mg$-type if there exists an element $g$ in $\mg$ such that $P_1.g = P_2$.  {\it A Delone set of  finite $\mg$-type} is a Delone set such that for each positive number $T$, there exists   only a finite number of $\mg$-types of patches in $X$ with diameter smaller than $T$.  The Lie group $\mg$ has a natural right action on  the set of Delone sets of finite $\mg$-type. Furthermore, any metric  on $\mg$  induces a metric  $\delta_\mg$ on the orbit $X.\mg$ of a given Delone set $X$ of finite $\mg$-type  defined as follows:

\noindent Fix a point $0$ in $\mg/\Gamma$ and let $B_\epsilon(0)$ stand for the open ball with radius $\epsilon$ and center $0$ in $\V$. Consider two Delone sets $X.g_1$ and  $ X.g_2$ in $X.\mg$. Let $A$ denote the set of $\epsilon \in ]0,1[$ such that there exist $g$ and 
$g'$ $\epsilon$-close to the identity $e$ in $\mg$ such that $X.g_1g\cap B_{1/\epsilon}(0) = X.g_2g'\cap 
B_{1/\epsilon}(0)$, then
$$ \delta_\mg (X.g_1, X.g_2)= \inf A \ \ \ \ {\rm if} \ \ A\neq \emptyset$$
$$\delta_\mg (X.g_1, X.g_2)=1  \ \ \ \ {\rm if} \ \ A = \emptyset \ .$$

\noindent Hence the diameter of  $X.\mg$ is bounded by $1$ and the $\mg$-action  on $X.\mg$ is continuous. The {\it continuous hull} $\Omega_\mg(X)$ of a Delone set of finite $\mg$-type $X$ is the completion of the metric space $(X.\mg, \delta_\mg)$. It is straightforward to check (see for instance \cite{KP}) that $\Omega_\mg(X)$ is a compact metric space and that any element in $\Omega_\mg(X)$ is a Delone set whose $\mg$-type of patches are those of $X$. Thus, the group $\mg$ acts on $\Omega_\mg(X)$  and  the dynamical system $(\Omega_\mg(X), \mg)$ possesses (by construction) a dense orbit (namely the orbit $X.\mg$). 

A Delone set of finite $\mg$-type $X$ is {\it repetitive} if for any patch $P$ in $X$ there exists a radius $\rho(P)>0$ such that any ball with radius $\rho(P)$ intersects $X$ in a patch that contains a patch with same $\mg$-type as $P$.
This last property can be interpreted in the dynamical system framework (see for instance \cite{KP}): 

\begin{prop}\label{min} The dynamical system $(\Omega_\mg(X), \mg)$ is minimal (i.e all its orbits are dense) iff the Delone set of finite type $X$ is repetitive.
\end{prop}
A Delone set $X$  is {\it strongly aperiodic} if there exist no $g\neq e$  $\mg$  and  no $Y$ in $\Omega_\mg(X)$ such that $Y.g =Y$.

\bigskip
\begin{remark} \label{re}{\rm  It is not clear  a priori that strongly aperiodic Delone sets of finite $\mg$-type exist for any homogenous space $\mg/\Gamma$ of any connected Lie group $\mg$. When the Lie group is the group of translations in $\mr^2$, the Penrose Delone set  is a standard example.   When  the Lie group is the group of   direct isometries $\me^2$ acting on $\mr^2$, the Penrose Delone set is no longer strongly aperiodic  since some elements in its Hull  have a five-fold symmetry. However strongly aperiodic examples can be made by a slight modification of the Penrose Delone set.   The Pinwheel Delone set \cite{Radin} is an example of repetitive  Delone set of finite $\mg$-type ; no Delone set in its hull is fixed by an infinite order element in $\mg$, however there are elements with a two-fold symmetry.    S. Mozes \cite{Mozes} has proved that when $\mg$ is either a non compact Lie group of rank at least 2, or  $Sp(n, 1)$ for $n\geq 2$ or ${\bf F}^{-20}_4$ there exist also  in $\mg/\Gamma$,  Delone sets of finite $\mg$-type which are not fixed by any element in $\mg$ of infinite order. Finally, very recently ,C. Goodman-Strauss \cite{gs1}, \cite{gs2}, has shown a similar result when the Lie group is $PSL(2, \mr)$ and the homogeneous space is the Poincar\'e disk.}
\end{remark}

\bigskip
Fix a point $0$ in $\mg/\Gamma$  such that $\Gamma$ fixes $0$. Consider the subset  $\Omega_{\mg}^0(X)$ in $\Omega_{\mg}(X)$ of the Delone sets $Y$ such that $0$ belongs to $Y$.  The group $\mg$ does not act on $\Omega_{\mg}^0(X)$ but the compact subgroup $\Gamma$  does.  The {\it canonical transversal} of a Delone set $X$ is the quotient space 
$\Omega_{\mg}^0(X)/ \Gamma$.  
 \begin{prop}The canonical transversal of a strongly  aperiodic  Delone set of finite $\mg$-type equipped with the quotient metric, is a totally disconnected metric space. Furthermore, if the Delone set is repetitive, it is  a metric Cantor set.
\end{prop}

\noindent{\it Proof:}  Two elements $Y_1$ and $Y_2$ in  $\Omega_{\mg}(X)$ are close one to the other if in a big  ball surrounding $0$, the two associated patches are images one of the other  by the action of an element in $\mg$ close to $e$. If $Y_1$ and $Y_2$ are in  $\Omega_{\mg}^0(X)$, this implies that the two patches are images one of the other  by the action of an element in $\Gamma$ close to $e$.  Since $X$  has finite $\mg$-type, we have constructed this way a basis of clopen neighborhoods of $\Omega_{\mg}^0(X)/ \Gamma$ which  consequently is totally disconnected. The repetitivity yields directly the  fact that  $\Omega_{\mg}^0(X)/ \Gamma$ is perfect. 

Two Delone sets $X$ and $Y$ are {\it $\mg$-equivalent (resp. $\mg$-orbit equivalent)} if the two dynamical systems $(\Omega_\mg(X), \mg)$ and $(\Omega_\mg(Y), \mg)$ are conjugate (resp. orbit equivalent) {\em  i.e.} there exists a homeomorphism $h:\Omega_\mg(X)\to \Omega_\mg(Y)$ that commutes with the $\mg$-action (resp. that maps $\mg$-orbits onto $\mg$-orbits).  

\bigskip
The following two theorems relate Delone sets and solenoids.
\begin{teo}
\label{sol0}
\noindent  
 Let $\mg$ be a connected Lie group, $\Gamma$ its maximal compact subgroup and let  $X$ be a strongly aperiodic Delone set of finite $\mg$-type in $\mg/\Gamma$. Then,  the continuous hull $\Omega_{\mg}(X)$   inherits a structure of an expansive  $\mg$-solenoid and the $\mg$-action on $\Omega_{\mg}(X)$ coincides with the canonical $\mg$-action of the $\mg$-solenoid. This action is minimal if and only the Delone set is repetitive. \end{teo}

\noindent {\it Proof: }  Consider a strongly aperiodic, repetitive,  Delone set $X$ in $\mg/\Gamma$ of finite $\mg$-type and choose a section $s:  \Omega_{\mg}^0(X)/ \Gamma\to   \Omega_{\mg}^0(X).$ Any point $X'$ in $\Omega_{\mg}(X)$ reads  in a unique way $X'= X".g $ where $X"$ is in $s(\Omega_{\mg}^0(X)/ \Gamma)$ and $g$ is an element in $\mg$. It follows that $\Omega_{\mg}(X)$ can be covered by a finite union of sets  $U_i=\phi_i(V_i\times T_i)$ where:
\begin{itemize}
\item  $T_i$ is a clopen set in $s(\Omega_{\mg}^0(X)/\Gamma)$;
\item  $V_i$ is an open set in $\mg$ which reads $V_i = \Gamma\times W_i,$  where $W_i$ is an open set in $\mg/\Gamma$;
\item  $\phi_i: V_i\times T_i\to \Omega_{\mg}(X)$ is the map defined by $\phi_i(v, t) = t.v^{-1}$. 
\end{itemize}
 Since there exists a finite partition of  $\Omega_{\mg}^0(X)/\Gamma$ (and of any of its clopen sets) in clopen sets with arbitrary small maximal  diameter, it is possible to choose this  diameter  small enough so that:
\begin{itemize}
\item the maps $\phi_i$ are homeomorphisms onto their images;
\item whenever $ Z\in U_i\cap U_j$, $Z=\phi_i(v, t)=\phi_j(v', t')$, the element $v^{-1}v'$ in $\mg$ is  independent of the choice of $Z$ in  $U_i\cap U_j$, we denote it by $g_{i, j}$.
\end{itemize}
The transition maps read $(v', t') = (v.g_{i, j}, s^{-1}(s(t) g_{i,j})$.
It follows that the boxes $U_i$ and charts $h_i=\phi_i^{-1}: U_i\to V_i\times T_i$ define a $\mg$-solenoid structure on $\Omega_{\mg}(X)$.  By construction, the $\mg$-action on this $\mg$-solenoid coincides with the $\mg$-action  on $\Omega_{\mg}(X)$. The $\mg$-orbits of $\Omega_{\mg}(X)$ are mapped to leaves of the $\mg$-solenoid structure. 

\noindent The $\mg$-action on $\Omega_{\mg}(X)$ is expansive. Indeed consider two Delone sets $Y_1$ and $Y_2$ in $s (\Omega_{\mg}^0(X)/ \Gamma)$ and within a distance $\epsilon$ small enough so that they coincide on a ball centered at $0$ and that their restriction to this ball $ B$  is fixed by no element in $\mg -\{e\}$. If $Y_1.g$ remains $\epsilon$-close to $Y_2.g$ when $g$ runs over a ball $B_\mg$, then $Y_1$ and $Y_2$ coincide on $B.B_\mg$. Thus,  if  the two $\mg$-orbits remain $\epsilon$-close,  $Y_1$ and $Y_2$ coincide. 

\noindent
The equivalence between the minimality of these two dynamical systems is plain and is equivalent to the repetitivity of the Delone set $X$ thanks to Proposition \ref{min}.  

\bigskip We say that  a  finite collection of verticals $V$ chosen in boxes of an atlas of  a $\mg$-solenoid $M$ is {\it well distributed} if there exists a constant $R>0$ such that, for any point $\tilde x$ on any leaf of $M$, there exists $\tilde y = x.g$ in $V$ where $\pi(g)$ is  $R$-close to $0$ in $\mg/\Gamma$  and that is has {\it no local symmetries} if furthermore, $x.g.\gamma$ is not in $V$ for any  $\gamma$ in $\Gamma$ different from $e$.  

\noindent We say that   $V$ is {\it of type $(\epsilon, R)$ } if:
\begin{itemize}
\item [(i)] $V$ is  well distributed with  constant $R>0$ and no local symmetries;
\item [(ii)] the diameters of the verticals in $V$ is smaller that $\epsilon$;
 \item  [(iii)] for any pair of verticals $\Delta_1$ and $\Delta_2$ in $V$ and any $\tilde x$ in  $\Delta_1$ and $g$ in $\mg$ such that $\tilde x g$ is in $\Delta_2$, then for any $\tilde y$ in $\Delta_1$, $\tilde y.g$ is in $\Delta_2$ or $\tilde y.g$ is not in $V$.
\end{itemize}
Since we have a lot of freedom in choosing the verticals, it is clear that finite collection of verticals   of type $(\epsilon, R)$ exists on each $\mg$-solenoid for each $R>0$ and $\epsilon >0$. 

\noindent When the $\mg$-solenoid $M$ is expansive, we say that a finite collection of verticals is {\it adapted to $M$ } if it is of type $(\epsilon, R)$ for some 
$R>0$ and $\epsilon >0$ and if furthermore, for any vertical $\Delta$ in $V$ 	and any pair of points $\tilde y_1$ and $\tilde y_2$ in $\Delta$,  the distance between $\tilde y_1.g$ and $\tilde y_2.g$ is smaller than the expansivity diameter $\epsilon (M)$ as long as $\pi(g)$ remains  $R$-close to $0$ in $\mg/\Gamma$. Here again, for any expansive $\mg$-solenoid, it is always possible to find an adapted collection of verticals.
 
\begin{teo}\label{sol}
Let $\mg$ be a connected Lie group, $\Gamma$ its maximal compact subgroup and $M$ a  $\mg$-solenoid. Then:
\begin{enumerate}
\item For any finite collection $V$  of relatively dense verticals chosen in boxes of an atlas of $M$,   the intersection of $V$ with any leaf  $L$ of $M$, defines (by identifying $L$ with  $\mg$ and projecting on  $\mg/ \Gamma$)  a  Delone set $Y$ of finite $\mg$-type on $\mg/\Gamma$. 
\item  If $M$ is expansive and $V$ is adapted to $M$,  $Y$ is strongly aperiodic.
\item  If $M$ is expansive and minimal and $V$ is adapted to $M$, the dynamical systems $(\Omega_\mg (Y), \mg)$ and $(M, \mg)$ are conjugate and thus all the Delone sets constructed in this way are $\mg$-equivalent and repetitive.
\end{enumerate}
\end{teo}

\noindent {\it Proof: } 

\noindent (1)- Let $\pi: \mg\to \mg/\Gamma$ be the standard projection. Let $V$ finite union of verticals in the $\mg$-solenoid and consider the projection $Y $ on $\mg/\Gamma$ of the intersection $V\cap L$ where $L$ is any leaf of the solenoid: $$  Y = \pi (V\cap L).$$  Whenever $V$ is well distributed, it is clear that $Y$ is a Delone set. Fix a number $N >0$ and assume that $\tilde y_0$ is in $V\cap L$ and denote $y_0 =\pi (\tilde y_0)$ its projection  in $Y$, the other intersection points in the ball with center $y_0$ and radius $N$ being  $y_1 =\pi (\tilde y_1) , \dots, y_n =\pi (\tilde y_n)$, where $\tilde y_1 =\tilde y_0.g_1,\dots , \tilde y_n =\tilde y_0.g_n$. Consider a point $\tilde z_0$  in  $V$ close enough to $\tilde y_0$.  Property $(5)$ above insures that, for $i= 1, \dots, n$,  $\tilde z_0.g_i$ is in $V$.  Thus, the set of points $z_0 =\pi (\tilde z_0),\dots z_n =\pi (\tilde z_n)$  is an isometric copy of the set $y_0, \dots, y_n$. By choosing a finite partition of $V$ in clopen sets with small enough diameters, we get that there exists only a finite number of $\mg$-type for patches with size smaller than $N$. Since this is true for all $N>0$, $Y$ is a Delone set of finite $\mg$-type. 

\noindent (2)- Assume now that $V$ is adapted to $M$.
We can construct the map $\Phi$ from the orbit $\tilde y_0.\mg$ in $(M, \mg)$ to the orbit $Y.\mg$  in $(\Omega_\mg(Y), \mg)$ that maps $\tilde y_0$ to $Y$ and makes the two $\mg$-actions commute.  The map $\Phi$ is continuous for the topologies induced by $M$ on $\tilde y_0.\mg$ and induced by $\Omega_{\mg}(Y)$ on $Y.\mg$. Let us show that $\Phi:\tilde y_0.\mg\to Y.\mg$  is a bijection. For this purpose we use the expansivity on $M$. Assume that $\Phi$ is not bijective, this means that there exists $g\neq$ in $\mg$  such that the Delone sets associated with $\tilde y_0$ and  $\tilde y_0.g$ coincide. Since $V$ is has no local symmetries, no element in $\Gamma$ different from $e$ let $Y$ globally invariant. Thus $g$ must be outside the maximal compact group.
  Notice that $\tilde y_0.g^n$  is also in $V$ for all $n>0$. Choose two points $\tilde a$ and $\tilde b$ in the sequence $\tilde y_0.g^n$ which are on a same vertical $\delta$ in $V$. We now use the fact that $V$ is adapted to $M$. This implies (point $(iii)$) that if $\tilde a.g'$ is in some vertical of $V$ for some $g'$ in $\mg$, $\tilde a.g'$ is in the same vertical and thus the $\mg$-orbits of $\tilde a$ and $\tilde b$ remain $\epsilon (M)$-close. The expansivity implies that  $\tilde a$ and $\tilde b$  coincide. A contradiction. It follows that the map $\Phi$ extends to a homeomorphism from the closure of $\tilde y_0.\mg$  in $M$ to $\Omega_{\mg}(Y)$. Consequently  $\mg$ has a free action  $\Omega_{\mg}(Y)$ and  $Y$ is strongly aperiodic.

\noindent (3)- 
If $M$ is minimal, the dynamical systems $(\Omega_\mg (Y), \mg)$ and $(M, \mg)$ are conjugate by $\Phi$ and thus $Y$ is repetitive. 
\section{Box decompositions and tower systems}
\subsection {Box decomposition }
 A {\it box decomposition} of a $\mg$-solenoid   $M$ is a finite collection of boxes $B_1, \dots B_n$ such that any two boxes are disjoint and the closure of the union of all boxes covers the whole solenoid. A box decomposition is {\it polyhedral} if each box $B_i$ reads in a chart $\pi^{-1}(U_i)\times T_i$ where $U_i$  an open convex geodesic polyhedron in $\mg/\Gamma$. Since the intersection of two polyhedral boxes is a finite union of polyhedral boxes, it is plain to check that every  $\mg$-solenoid admits a polyhedral box decomposition. In general, box decompositions associated with a $\mg$-solenoid   $M$ are not unique; however if we fix   a finite collection of verticals $V$ in the solenoid, there exists a canonical polyhedral box decomposition associated with the pair  $(M, V)$. In particular, given strongly aperiodic , repetitive, Delone set $X$ of $\mg$-finite type,  there exists a canonical polyhedral box decomposition of its hull $\Omega_{\mg}(X)$. We call this particular box decomposition, the {\it Voronoi box decomposition} and define it as follows:

\noindent Let $V$ be a well distributed finite union of verticals in a  $\mg$-solenoid $M$,  and $Y$ the projection  on $\mg/\Gamma$ of the intersection $V\cap L$ where $L$ is any leaf of the solenoid. We know from Theorem \ref{sol} that $Y$ is  a Delone set of $\mg$-finite type in the leaf $\mg/\Gamma$.  Consider the Voronoi tiling associated with $Y$ in $L$. {\em  i.e.} the polyhedral tiling   defined as follows:

\noindent For any $y$ in $Y$ the {\it Voronoi cell} $V_y$ of $y$ is  the open set of points in  $\mg/\Gamma$  which are closer to $y$ than to any other points in $Y$. The family $\{ \overline{V_y}\,; \; y\in Y\}$
defines a polyhedral tiling of  $\mg/\Gamma$ , in which the tiles are convex and meet face-to-face. Again from Theorem \ref{sol},  there exists a finite partition of $V$ in clopen sets $C_1, \dots, C_n$, such that for any pair of points $\tilde y_1$ and $\tilde y_2$ in a same $C_i,$  $V_{\pi(\tilde y_2)}$ is obtained from $V_{\pi(\tilde y_1)}$  by the action of an element in $\mg$. The boxes of the Voronoi box decomposition read in charts $\pi^{-1} (V_y)\times C_i$ where $y$ is in $Y$ and $i=1, \dots, n$.

\subsection{Tower system}
The {\it vertical boundary} of a polyhedral box which read in a chart $\pi^{-1}(P)\times T$ where $P$ is a polyhedron in $\mg/\Gamma$ and $T$ a totally disconnected set, reads in the same chart $\pi^{-1}(\partial P)\times T$, where $\partial P$ stands for the boundary of $P$.  A polyhedral box decomposition $\B_2$ is {\it zoomed out} of a polyhedral box decomposition $\B_1$ if:
\begin{itemize}
\item [(1)] for each point $x$ in a box $B_1$ in $\B_1$ and in  a box $B_2$ in $ \B_2$, the vertical of $x$ in $B_2$ is contained in the vertical of $x$ in $B_1$;
\item [(2)] the vertical boundaries of the boxes of $\B_2$ are contained in the vertical boundaries of the boxes of $\B_1$;
\item [(3)] for each box $B_2$ in $\B_2$, there exists a box $B_1$ in $\B_1$ such that $B_1\cap B_2\neq \emptyset $ and  the vertical boundary of $B_1$ does not intersect the vertical boundary of $B_2$. 
\item [(4)]  if a vertical in the vertical boundary of  a box in $\B_1$   contains a point in a vertical boundary of a box in $\B_2$, then it contains the whole vertical\footnote{In substitution theory, this condition is called {\it forcing the border} \cite{KP}. In the study of Williams attractors it is called {\it flattening condition} \cite{W2}. }.
\end{itemize}
A {\it tower system} of a  $\mg$-solenoid is a sequence of polyhedral  box decompositions $(\B_n)_{n\geq 1}$ such that for each $n\geq 1$, $\B_{n+1}$ is zoomed out of  $\B_n$.
\begin{teo}\label{tour} Let $\mg$ be a connected Lie group, then any  $\mg$-solenoid admits a tower system.\end{teo}

\noindent{\it Proof:} Consider a triplet $(M, V_1, x)$ where $M$ is a  $\mg$-solenoid,   $V_1$ is  a well distributed finite collection of verticals and $x$ a point in $V_1$. 

\noindent{\it Step 1:} Fix $\delta_1 =1$. The first box decomposition $\B_1$ is the Voronoi box decomposition associated with $(M, V_1)$.

\noindent{\it Step 2:} Consider the box $B_{1, 1}$ in $\B_1$ that contains $x$ and choose in  $B_{1, 1}\cap V_1$ a clopen set $V_2$ which contains $x$ and has diameter smaller than $\delta_2$.  Consider the Voronoi box decomposition $\B_2'$ associated with  $(M, V_2)$. The box decomposition $\B_2'$ is not zoomed out of $\B_1$. Indeed, it satisfies point 1 of the definition, and if $\delta_2$ is chosen small enough it satisfies also points 3 and  4, but it certainly does not satisfy point 2. To construct a box decomposition which satisfies also point 2   we have to introduce a slight modification in the construction of the Voronoi box decomposition. More precisely, let $Y_2= \pi(V_2\cap L)$ where $L$ is a leaf of $M$.   Instead of considering  the Voronoi tiling associated with $Y_2$ in $L$ and for each point $y$ in $Y_2$ its Voronoi cell $V_y$, we construct the best approximation possible\footnote{The mild ambiguity inherent with this construction will turn to be irrelevant because of the $\mg$-finite type hypothesis.} of the Voronoi tiling by changing each Voronoi cell $V_y$ by a reunion $\tilde V_y$ of cells of the Voronoi tiling associated with the Delone set $ Y_1=\pi(V_1\cap L)$. The new boxes of the decomposition are then constructed following the same rules.

\noindent {\it Step 3:} We iterate this construction by  considering the box $B_{1, 2}$ in $\B_2$ that contains $x$ and choose in  $B_{1, 2}\cap V_2$ a clopen set $V_3$ which contains $x$ and has diameter smaller than $\delta_3$ etc.

\section {$\mg$-branched manifolds}
In this section we keep the same notations:  $\mg$ is a connected Lie group, $\Gamma$ its maximal compact subgroup. The quotient space    $\mg/\Gamma$ (right classes)  is equipped with a Riemannian metric such that $\mg$ has a transitive right action  by isometries on $\mg/\Gamma$ and $\pi:\mg\to \mg/\Gamma$ is the standard projection. 

\subsection{Local models}

For $r> 0$ and $x$ in $\mg/\Gamma$,  we  consider the open ball with radius $r$ centered at $x$, $B_r(x)$. The {\it  $\mg$-ball of type 1} $B( x, r)$ is the open set  $\pi^{-1}(B_r(x))$. A {polyhedral decomposition} of $B_r(x)$ is a finite collection of polyhedral open cones centered at $x,\,$ $\,({\cal C}= \{C_1, \dots , C_n\})$ which are pairwise disjoint and whose closures cover the ball $B_r(x)$.
For $p\geq 1$, a {\it  $\mg$-ball of type $p$}  is defined in several steps:
\begin{itemize}
\item Consider in the ball $ B_r(x)$ a finite collection of polyhedral decompositions $\C_1, \dots, \C_p$ such that
if two cones in the same $\C_i$ have a $(d-1)$-face in common, then there exits $j\neq i$ such that one of the two cones appears also in  $\C_j$ and the other one does not. 
\item Define in the disjoint union $\sqcup_{i=1}^{i=p} B_r(x)$ the equivalence relation:
$(y, i) \sim (y', j)$ iff $\C_i$ and $\C_j$ share a same cone whose closure contains $y=y'$.
\item In $\sqcup_{i=1}^{i=p}\pi^{-1}( B_r(x))$,  consider the equivalence relation $\approx$ defined by $\tilde y \approx \tilde y'$ if and only if $\pi(\tilde y) \sim \pi (\tilde y')$.  
\item The {\it  $\mg$-ball of type $p$}, $B(x, r, \C_1,\dots, \C_p)$,  is the set $\sqcup_{i=1}^{i=p}\pi^{-1}( B_r(x))/\approx$.\end{itemize}
The  $\mg$-balls, as defined above, are said {\it centered at $x$ and with radius $r$}. Notice that the projection $\pi$ induces a projection from $\sqcup_{i=1}^{i=p}\pi^{-1}( B_r(x))/\approx$ to   $\sqcup_{i=1}^{i=p} B_r(x)/\sim$ that we still denote $\pi$.  Let  $\pi_1:  \sqcup_{i=1}^{i=p} B_r(x)\to  \sqcup_{i=1}^{i=p} B_r(x)/\sim$  be the canonical projection. This quotient space can be seen as a collection of $p$ copies of the ball $B_r(x)$ glued one with the other along polyhedral cones centered at $x$ and is stratified as follows: 

\noindent For $\tilde y$ in $\sqcup_{i=1}^{i=p} B_r(x)/\sim$  and $y$ in $\pi_1^{-1}(\tilde y)$, let $\nu(y)$ be the number of cones whose closure contains $y$ and $\nu(\tilde y)$ be the maximum of all the $\nu(y)$ when $y$ runs over $\pi_1^{-1}(\tilde y)$.  For  $l$ in $\{1, \dots, g \}$, where $g = dim \mg$,  $ {\cal V}_l $  is the set of points $\tilde y$ in $\sqcup_{i=1}^{i=p} B_r(x)/\sim$ such that $\nu(\tilde y)=  g-l+1$ and ${\cal V}_0$ is the set of points $\tilde y$ in $\sqcup_{i=1}^{i=p} B_r(x)/\sim$ such that $\nu(\tilde y) \geq  g+1$.  By denoting $cl(A)$ the closure of a set $A$, we have:
\begin{itemize}
\item $\sqcup_{i=1}^{i=p} B_r(x)/\sim\,= \, cl({\cal V}_g) \, =\,\cup_{l=0}^{l=g} {\cal V}_l$;
\item  for each $l$ in $\{1, \dots, g\}$, $\,cl({\cal V}_l)\, =\, {\cal V}_l\cup {\cal V}_{l-1}$ and ${\cal V}_l\cap {\cal V}_{l-1}=\emptyset$;
\item $cl({\cal V}_0) = {\cal V}_0$ contains the point $x$.
\end{itemize}

This stratification yields  a stratification of the $\mg$-ball. The {\it singular locus} of the $\mg$-ball with type $p$ is the set $\pi^{-1}(cl(\V_{g-1}))$. For  $l$ in $\{0, \dots, g\}$, $\pi^{-1}(\V_l)$ is called the {\it $l$-stratum} of the $\mg$-ball. Any $\mg$ ball inherits a smooth structure, a right $\Gamma$-action, a local $\mg$-action defined outside the singular locus and a Riemannian metric which is invariant under the local $\mg$-action where it is defined.

\noindent Let $proj$ be the standard projection $proj: B( x,r,  \C_1,\dots, \C_p)\to B(x, r)$. A {\it $\mg$-sheet} in the  $\mg$-ball $B( x,r,  \C_1,\dots, \C_p)$ is the image of $B(x, r)$  by a section of $proj$.

\noindent The following lemma, whose proof is straightforward, describes the situation around the center of a $\mg$-ball of type $p$.
\begin{lem} In  a $\mg$-ball with type $p$, there exists a neighborhood of each point $y$ which is a   $\mg$-ball centered at $ y$ of type $q$ where $q$ in $\{1, \dots , p\}$.  Furthermore, if $y$ is not in the singular locus,  then $q=1$.
\end{lem}
 A map $\tau: B( x_2, r_2,  \C^2_1,\dots, \C^2_{p_2})\to B( x_1, r_1, \C^1_1,\dots, \C^1_{p_1})$ is a { \it $\mg$-local submersion} if:
\begin{itemize} 
\item The map $\tau$ preserves the smooth structure {\em  i.e.} it maps:
\begin{itemize}
\item the singular locus  into the singular locus;
\item the  $\mg$-sheets into $\mg$-sheets.
\end{itemize}
\item There exists $g$ in $\mg$ such that  $\tau$ projects on the right multiplication by  $g$ on $B(x_2, r_2)$ {\em  i.e.} $proj\circ\tau(.) =proj(.).g$.
\end{itemize}
A $\mg$-local submersion is a {\it $\mg$-local isometry} if it is one-to-one. It this case, the two $\mg$-balls have the same radius and the same type and the inverse map   is also  a  $\mg$-local isometry.

\subsection{Global models}
Consider a compact metric space $S$ and assume there exist a cover of $S$ by open sets $U_i$  and  chart homeomorphisms $h_i:U_i\to V_i$ where $V_i$ is an open set in some  $\mg$-ball. These open sets and homeomorphisms define an atlas of a {\it $\mg$-branched manifold} structure on $S$ if the transition maps  $h_{i,j}\, =\, h_j\circ h_i^{-1}$ satisfy the following property:

\noindent  For each point $y$ in $h_i(U_i\cap U_j)$,  there exists a neighborhood $B_y$ of  $y$  which is a $\mg$-ball contained in $h_i(U_i\cap U_j)$ and a neighborhood $B_{h_{i,j}(y)}$ of  $h_{i,j}(y)$  which is a $\mg$-ball contained in $h_j(U_i\cap U_j)$  such that the restriction of $h_{i,j}$ to $B_y$ is a local $\mg$-isometry from $B_y$ to  $B_{h_{i,j}(y)}$. 

\noindent Two atlas are {\it equivalent} if their union is again an atlas. A {\it $\mg$-branched manifold}
is the data of a metric compact space $S$ together with an equivalence class of atlas $\mathcal{A}$.

\begin{lem} \label{gamma}
A $\mg$-branched manifold is canonically equipped with a smooth structure,  a Riemannian metric, a stratification which is $\Gamma$-left invariant, and a local right $\mg$-action defined outside of the singular locus which preserves the Riemannian metric. 
\end{lem}
{\it Proof:} 

\noindent The smooth curve on a $\mg$-branched manifold is a curve that, when read in a chart,  remains in a $\mg$ -sheet and is mapped by $proj$ onto a smooth curve in some $\mg$. 

\noindent In a chart, the Riemannian metric is the pullback of standard invariant metric on $\mg$ by $proj$. It is clear that this choice is coherent with respect to changes of charts.

\noindent The stratification $\S$ is defined as follows: 
\noindent For  $l$ in $\{0, \dots, g\}$, the {\it $l$-stratum of $\S$} is the set of points which, when read in a chart,  have a neighborhood which is a $\mg$-ball and belong to the $l$-stratum of this $\mg$-ball.  This property is independent on the chart and is $\Gamma$-left invariant. A {\it $l$-region} is a connected component of the $l$-stratum. The finite partition of $S$ in $l$-regions, for $l$ in $\{0, \dots, g\}$, is called the {natural stratification of the $\mg$-branched manifold}. The union of all the $l$ -regions for $l$ in $\{0, \dots, d-1\}$  forms the {\it singular locus $Sing(S)$ of $S$}. 

\noindent Outside of the singular locus, a neighborhood of a point is modeled on a neighborhood in $\mg$, it is thus possible to define a local $\mg$ action. 

\bigskip 
Let $S$ be a $\mg$-branched manifold and $x$ a point in $S$.    The 
{\it injectivity radius of $x$}  is the supremum of the radii of all $\mg$-sheets centered at $x$.  The {\it injectivity radius of $S$}, denoted $inj(S)$,  is the  infimum  over all $x$ in $S$ of the injectivity radius of $x$. The {\it $\mg$-action radius} is the infimum of all connected components of $S\setminus Sing(S)$ of the supremum of the radius of a $\mg$-ball embedded in the connected component, we denote it $size (S)$. 

Let $S_1$ and $S_2$ be two $\mg$-branched manifolds.  A continuous surjection $\tau:S_2\to S_1$ is a{ $\mg$-submersion} if for any $y$ in $S_2$  there exist a neighborhood $B_y$ of  $y$  which is a $\mg$-ball (when read in a chart) and a neighborhood $B_{\tau(y)}$ of  $\tau(y)$  which is a $\mg$-ball (when read in a chart)  such that the restriction of $\tau$ to $B_y$ is a $\mg$-local submersion from $B_y$ to  $B_{\tau(y)}$.  It is a {\it  strong $\mg$-submersion} if $\tau(B_y)$ is in a $\mg$-sheet of $B_{\tau(y)}$. 
\begin{lem}\label{inj}
 Let $S_1$ and $S_2$ be two $\mg$-branched manifolds and  $\tau:S_2\to S_1$ a $\mg$-submersion. Then $inj(S_2) \geq inj(S_1)$ and  $size (S_2) \geq size(S_1)$.
\end{lem}

\noindent{\it Proof:} A $\mg$-submersion $\tau: S_2\to  S_1$ preserves the metric of the $\mg$-branched surfaces. It follows that for any $\mg$-sheet $\F_1$ in $S_1$, there exists a $\mg$-sheet $\F_1$ in $S_2$ (isometric to $\F_1$) that is mapped on $\F_1$ by $\tau$, which implies the lemma.

\subsection{Projective limits}
For $n\geq 1$, let $S_n$ be a sequence of $\mg$-branched manifolds and $\tau_n : S_{n+1}\to S_n$ a sequence of strong $\mg$-submersions. Let us recall that the elements of the {\it projective limit } $\lim_\leftarrow (S_n, \tau_n)$
consist in the elements  $(x_1, x_2, \dots, x_n, \dots)$ in the product $\prod_{n\geq 1} S_n$ such
that $\tau_n(x_{n+1}) = x_n$ for all $n\geq 1$. Since all the $S_n$'s are compact,   $\prod_{n\geq 1} S_n$ equipped with the product topology is compact. The projective limit being a closed subset of this product, is also compact. 
For every $i\geq 1$, we denote by $p_i$ the  natural   continuous map
$p_i: \lim_\leftarrow (S_n, \tau_n) \to S_i$.

 \begin{prop} For $n\geq 1$, let $S_n$ be a sequence of $\mg$-branched manifolds and $\tau_n : S_{n+1}\to S_n$ a sequence of strong $\mg$-submersions. Then, the projective limit   $\lim_\leftarrow (S_n, \tau_n)$ is equipped with  a natural local $\mg$-action.
\end{prop}

\noindent{\it Proof:} As we have already observed (Lemma \ref{gamma}), there exists a local $\mg$-action on each $\mg$-branched manifold $S_n$ outside of the singular locus. It is straightforward to observe that this action commutes with the $\mg$-submersions $\tau_n$. Consider now a point $(x_1, x_2, \dots, x_n, \dots)$   in the projective limit  
$\lim_\leftarrow (S_n, \tau_n) $.   From Lemma \ref{inj}, there exist a positive number $r>0$ and a sequence of $\mg$-sheets $\F_n$ in $S_n$ centered at $x_n$, with radius $r$ such that $\tau_n(\F_{n+1}) = \F_n$ for all $n \geq 1$. When the maps $\tau_n$ are $\mg$-submersion, this sequence $\S_n$ is not necessarily unique. However, it is unique when the maps $\tau_n$ are required to be strong $\mg$-submersions.  For $g$ in $\Gamma\times B_r(0)\subset \mg$, and for each $n\geq 1$, we define in a chart $ x_n.g = e_n(g)$ where $e_n: \Gamma\times B_r(0)\to S_n$ is the isometric embedding defining $\F_n$. It is plain to check that   $\tau_n(x_{n+1}.g) = x_n.g$ for all $n\geq 1$.  
\subsection {$\mg$-solenoids and projective limits of $\mg$-branched manifolds}
The following theorem relates  $\mg$-solenoids  and projective limits of $\mg$-branched manifolds.
\begin{teo}\label{branched} Let $\mg$ be a connected Lie group, $\Gamma$ its maximal compact subgroup  and $M$  a  $\mg$-solenoid:
\begin{enumerate}
\item there exist a sequence $S_n$  of $\mg$-branched manifolds and a sequence  of strong $\mg$-submersions $\tau_n : S_{n+1}\to S_n$ such that $M$ is homeomorphic to the projective limit  $\lim_\leftarrow (S_n, \tau_n)$;
\item  the local $\mg$-action on $\lim_\leftarrow (S_n, \tau_n)$ extends to a $\mg$-action and the homeomorphism  realizes a conjugacy between the two dynamical systems $(M, \mg)$  and $(\lim_\leftarrow (S_n, \tau_n), \mg)$;
\item $\lim_{n\to +\infty}inj(S_n) =+\infty$ and $\lim_{n\to +\infty}size(S_n) =+\infty$.
\item and all the $g$-regions of $S_n$ have the homotopy type of $\Gamma$; 
\item If $M$ is minimal then all the $S_n$'s are connected.
\end{enumerate}\end{teo}

\noindent{\it Proof:} Consider a polyhedral box decomposition $\B_1$ of $M$ together with the following equivalence relation $\simeq$ that identify two points that are in a same vertical in the closure of a box of the decomposition. The whole construction we made in the previous section has been adapted so that the quotient space $M/\approx$ inherits a natural structure of $\mg$-branched manifold $S_1$. Let us denote by $p_1: M\to S_1$, the standard projection.  From Theorem \ref{tour}  any  minimal $\mg$-solenoid admits a tower system {\em  i.e.} a sequence of polyhedral  box decompositions $(\B_n)_{n\geq 1}$ such that for each $n\geq 1$, $\B_{n+1}$ is zoomed out of $\B_n$. We associate with this sequence the corresponding sequence of $\mg$-branched manifolds $S_n$ and maps $p_n: \Omega_{\mg}(X)\to S_n$. Remark that two points which are on a same vertical in a closed box in  $\B_{n+1}$ are also on a same vertical of a box in $\B_n$. This defines a map $\tau_n:  S_{n+1}\to  S_n$. It is clear that $\tau_n$ is a $\mg$-submersion and that $\tau_n\circ \pi_{n+1} =\pi_n$, for all $n\geq 1$. Furthermore, since each $\B_ {n+1}$ is zoomed out of $\B_n$, we know that a vertical in the vertical boundary of  a box in $\B_{n}$ contains or is disjoint from any vertical boundary of a box in $\B_{n+1}$. It follows that $\tau_n$ maps a small neighborhood of a point in the singular set of $\S_{n+1}$ into a $\mg$-sheet of $\S_{n}$ and thus it is a strong $\mg$-submersion.

\noindent At this point, we have constructed  a sequence $S_n$  of $\mg$-branched manifolds and a sequence  of strong $\mg$-submersions $\tau_n : S_{n+1}\to S_n$ and a map $p:\Omega_{\mg}(X)\to \lim_\leftarrow (S_n, \tau_n)$ defined by: $$p(x) =(p_1(x_1), \dots, p_n(x_n), \dots).$$ Its is plain to check that  $p$ is a homeomorphism and that it maps the local $\mg$ actions on onto the other. Since the $\mg$ action on $M$  is global, so is the local action on 
$\lim_\leftarrow (S_n, \tau_n)$. Thus, the  homeomorphism is a conjugacy between the two dynamical systems  $(M, \mg)$  and $(\lim_\leftarrow (S_n, \tau_n), \mg)$.

\noindent Since each $g$-region in $S_{n+1}$ is made of isometric copies of more than one region $S_n$, we immediately get that $\lim_{n\to +\infty}inj(S_n) =+\infty$ and $\lim_{n\to +\infty}size(S_n) =+\infty$. The fact that the box decompositions $\B_n$ are polyhedral yields that  the $g$-regions of $S_{\B_n}$ have the homotopy type of $\Gamma$.  The  fact that minimality implies connectedness is plain.

\section{ Transverse invariant measures of $\mg$-solenoid}
\subsection{Transverse invariant measure and Ruelle-Sullivan current}
Consider an atlas of a $\mg$-solenoid $M$ given by the charts $h_i:U_i\to V_i\times T_i$ where the $T_i$'s are totally disconnected metric sets and the $V_i$'s are open subsets in   $\mg$. { A \it finite transverse invariant measure} on $M$ (see \cite {Ghys}) is the data of a finite positive measure on each set $T_i$ in such a way that if $B$ is a 
Borelian set in some $T_i$ which is contained in the definition set of the transition map $\gamma_{ij}$ then:
$$\mu_i(B)\, =\, \mu_j(\gamma_{ij}(B)).$$

\noindent It is clear that the data of a transverse invariant measure for a given atlas provides another invariant measure for any equivalent atlas and thus gives an invariant measure on each vertical. Thus it makes sense to consider a transverse invariant measure $\mu^t$ of a $\mg$ solenoid. The fact that the leaves of a lamination carry a structure of $g$-dimensional manifold, where $g$ is the dimension of $\mg$ allows us to consider differential forms on $M$. 
A $k$-differential form on $M$  is the data of $k$-differential forms on the open sets $V_i$ that are mapped one onto the other by the differential of the transition maps $g_{ij}$. We denote by $A^k(M)$ the set of $k$-differential forms on $M$. A {\it foliated cycle} is a linear form from $A^g(M)$ to $\mr$ which is positive on positive forms and vanishes on exact forms.

\noindent There exists a simple way to associated with a transverse invariant measure  a foliated cycle. Consider a $g$-differential form $\omega$ in $A^g(M)$ and assume for the time being, that the support of $\omega$ is included in one of the $U_i$'s. In this case, the form can be seen as a form in $V_i\times T_i$. By integrating $\omega$ on the slices $V_i\times \{t\}$ we get a real valued map on $T_i$ that we can integrate against the transverse measure $\mu_i$ to get a real number $\C_{\mu^t}(\omega)$. When the support of $\omega$ is not in one of the $U_i$'s, we choose 
a partition of the unity $\{\phi_i\}_i$ associated with the cover of $M$ by the open sets $U_i$ and define:
$$\C_{\mu^t}(\omega)\,=\, \sum_i \C_{\mu^t}(\phi_i\omega).$$
It is clear that we have defined this way a linear form $ \C_{\mu^t} : A^g(M)\to \mr$ which does not depend on the choice of the atlas in $\L$ and of the partition of the unity. It is also easy to check that this linear form is positive for positive forms. The fact that $\C_{\mu^t}$ vanishes on closed form is a simple consequence of the invariance property of the transverse measure. The foliated cycle $\C_{\mu^t}$ is called {\it the Ruelle-Sullivan current} associated with the transverse invariant measure $\mu^t$. It turns out that the existence of a foliated cycle implies the existence of a transverse invariant measure (see \cite{Sullivan}) and thus both points of 
view: transverse invariant measure and foliated cycle are equivalent. We denote ${\cal M}^t(M)$ the convex set of transverse invariant measure on $M$.

\subsection{Transverse  and $\mg$-invariant measures}
Let us now use the fact  a  $\mg$-solenoid carries a $\mg$ action and assume that there exist a finite measure $\mu$ on $M$ which is invariant for the $\mg$-action. This invariant measure defines a transverse invariant measure. For any Borelian subset of a transverse set $T_i$:
$$\mu_i(B)\, =\, \frac{\mu (h_i^{-1}( V_i\times B))}{\lambda_\mg(V_i)},$$ 
where $\lambda_\mg$ stands for the  measure in $\mg$ which is invariant under the right action of $\mg$.  Conversely,
consider a transverse invariant measure $\mu^t$ of the  $\mg$-solenoid
$M$. Let $f: M\to \mr$ be  a continuous function  and assume for the time being that the support of $f$ is included in one of the $U_i$'s. In this case, the map $f\circ h_i^{-1}$ is defined on $V_i\times T_i$. By integrating $f\circ h_i^{-1}$ on the sheets $V_i\times \{t\}$ against the  measure $\lambda_\mg$ of $\mg$, we get a real valued map on $T_i$ that we can integrate against the transverse measure $\mu_i$ to get a real number $\int f d\mu$. When the support of $f$ is not in one of the $U_i$'s, we choose  a partition of the unity $\{\phi_i\}_i$ associated with the cover of $M$ by the open sets $U_i$ and define:
$$\int f d\mu\,=\, \sum_i \int f\phi_i d\mu.$$
It is clear that we have defined this way a finite  measure on $M$ which does not depend on the choice of the atlas in its equivalence class and of the partition of the unity and is invariant under the $\mg$-action.

\noindent Thus, for a $\mg$-solenoid  $M$, the following 3 points of view are equivalent:
\begin{itemize}
\item A finite transverse invariant measure;
\item a foliated cycle;
\item a finite $\mg$-invariant measure on $M$.
\end{itemize}
\subsection{Transverse invariant measures and homology}
Let us  first recall few elementary facts about the (cellular) homology of $S$. Consider the natural stratification  of $S$, ${\cal V}_0, \dots {\cal V}_l$, where, for $i=0, \dots g$, ${\cal V}_i$ is decomposed into a finite number $a_i$ of $i$-regions that we orient and order in an arbitrary way but for the orientation of the $g$-regions which is the natural orientation induced by the $\mg$-manifold structure. We refine this stratification by regions to get a cell decomposition and we orient each $g$-cell according to the orientation of the $g$-region that contains it. We  denote by $C_i(S,\mr)$  the free
$\mr$-module which has as (ordered) basis the set of ordered and oriented $i$-cells. By convention, for any oriented $i$-cell $e$, $-e = -1e$, and it consists of the same cell with the opposite orientation. 
 Notice that with the exception of the $g$-cells orientation, there are no canonical choices for the above orders and orientations, but these choices will be essentially immaterial in our discussion.

We define the  linear {\it boundary operator}:
$$\partial_{i+1} : C_{i+1}(S,\mr)\to  C_i(S,\mr)\ $$
which assigns to any $i+1$-region, the sum of the $i$-regions that are in its closure pondered with a positive sign (resp. negative) if the induced orientation fits (resp. does not fit) with the orientation chosen for these $i$-regions. It is clear that $\partial_i\circ\partial_{i+1}=0$. 

\noindent The  space $Z_i(S,\mr)= {\rm Ker}\ \partial_i$ is called the space of  $i$-{cycles} of $S$ and the space  $B_i(S,\mr)= \partial_{i+1}(C_{i+1}(S,\mr))$ is called the space of $i$-{boundaries} of $S$. In fact, $B_i(S,\mr)\subset Z_i(S,\mr)$ and $H_i(S,\mr)= Z_i(S,\mr)/ B_i(S,\mr)$ is the $i^{th}$ homology group of $S$. Notice that $H_g(S,\mr)=Z_g(S,\mr)$. 

A standard result of algebraic topology insures that (up to $\mr$-module isomorphism) $H_i(S,\mr)$ is a topological invariant of $S$ that coincides with the $i^{th}$ singular homology of $S$ (see for example \cite{Spa}). 

It is important to observe  that a $g$-chain $z$ is  a $g$-cycle ({\em  i.e.} $\partial_d(z)=0$) if and only if:
\begin{itemize}
\item two $g$-cells which are in a same region $R_i$ appear in $z$ with the same coefficient;
\item  the coefficients of $z$ satisfy the  ``switching rules''(or Kirchoff-like laws): this means that along every $g-1$-region $e$ of $S$ the sum of the weights on the germs of $d$-cells along $e$ on one side equal the sum of the
weights of the germs of $d$-cells on the other side.
\end{itemize}
Identifying the free $\mr$-module  which has as (ordered) basis the set of ordered and oriented $g$-regions  with of $\mr^{a_g}$,
 the homology group $H_g(S, \mr)$ is then the subspace  of $\mr^{a_g}$ defined by the switching rules.

 We denote by $H_g^+(S, \mr)$  the intersection of $ H_g(S, \mr)$ with the positive cone in $\mr^{a_i}$. 

\begin{teo} \label{measure}Let $M$ be a  $\mg$-solenoid  and  $S_n,\,$ $n\geq 1,$ be a sequence of $\mg$-branched manifolds and $\tau_n : S_{n+1}\to S_n$ a sequence of strong $\mg$-submersions such that  the two dynamical systems $(M,  \mg)$  and $(\lim_\leftarrow (S_n, \tau_n), \mg)$ are conjugate. Then:

$${\cal M}^t(M)\,\,\, \, \, {\rm is\,\, isomorphic\,\,  to}\,\,\,\, \lim_\leftarrow (H_g^+(S_n,\mr),  \tau_{n,\star})$$ 
where  $\tau_{n, \star}: H_g(S_{n+1}, \mr)\to H_g(S_n, \mr)$ is the map induced by $\tau_n$ on homology. 
\end{teo}

\noindent{\it Proof:} Let us choose one of these branched manifolds $S_n$ together with the natural projection $p_n: M\to S_n$. Let $R_{1, n},\dots, R_{p(n), n}$ be the ordered sequence of $g$-regions of $S_{n}$ equipped with the natural orientation and choose a $g$-region $R_{i,n}$ of $S_n$. The set $p_n^{-1} R$  reads in some chart $V_{i,n}\times T_{i, n}$. A transverse measure $\mu^t$ associates a weight  $w(R_{i, n}) =\mu^t(T_{i, n})$ with $R_{i, n}$ and thus defines an element $i_n(\mu_t)$ in $C_{g}(S_n,\mr)$. The fact that the transverse measure is invariant implies that the switching rules are satisfied, thus $i(\mu_t)$ is  in $ H_g(S_n, \mr).$ Since  a transverse invariant measures associate a positive weight to each  region, we conclude that $i(\mu_t)$ is  in  $ H_g^+(S_n, \mr)$. 

Let $R_{1, n+1}, \dots, R_{p(n+1), n+1}$ be the ordered sequence of  regions of $S_{n+1}$ equipped with the natural orientation and $\mr^p(n+1)$ the associated $\mr$-module.  By identifying $H_g(S_{n+1}, \mr)$ with a the subspace of $\mr^{p(n+1)}$ defined by the switching rules and  $H_g(S_{n}, \mr)$ with a the subspace of $\mr^{p(n)}$ defined by the switching rules,   we assign to  the linear map $$f_{n,*}:H_g(S_{n+1}, \mr)\to H_g(S_n, \mr),$$  a $n(p)\times n(p+1)$ matrix $A_n$ with integer non negative coefficients. The coefficient $a_{i, j, n}$   of the $i^{th}$ line and the $j^{th}$ column is exactly the number of pre images in $R'_j$ of a point in $R_i$.  Thus we have the relations:
$$w(R_{i, n})\, =\, \sum_{j=1}^{j=p(n+1)}a_{i, j, n}w(R_{j, n+1}).$$ 
for all $i=1,\dots, p(n)$ and all $j=1, \dots p(n+1),$ which exactly means that $\tau_{n,*}\circ i_n = i_{n+1}.$
Thus any finite, transverse, invariant measure can be seen as an element in $\lim_\leftarrow (H_g^+(S_n,\mr),  \tau_{n,\star})$.

\noindent Conversely, let $T_{1, 1}, \dots T_{p(1), 1}$ be the vertical associated with each $g$-region of the first $\mg$-branched manifold $S_1$. Since the $T_{i, 1}$'s are totally disconnected metric spaces, they can be cover by a partition in clopen sets with arbitrarily small diameters. Such a partition $\P$ is {\it finer} than another partition $\P'$ if the defining clopen sets of the first one are included in clopen sets of the second one. Consider a sequence of partitions $\P_n$, $n\geq 0$ of $T_{i,1}$ such that, for all $n\geq 0$ , $\P_{n+1}$ is finer than $\P_n$ and the diameter of the defining clopen sets of $\P_n$ goes to zero as $n$ goes to $+\infty$. A finite measure on $T_{i,1}$ is given by the countable data of non negative numbers associated with each defining clopen sets of each partition $\P_n$ which satisfy the obvious additivity relation. An element $y = (y_1, \dots, y_n, \dots)$ in $\lim_\leftarrow (H_g^+(S_n,\mr),  \tau_{n,\star})$ provides us such a sequence of partitions $\P_n$ defined as follows:
\begin{itemize}
\item let $\hat \pi_1:p^{-1}(R_{i, 1}) \to T_{i, 1}$ the chart map composed with the projection;
\item the clopen sets of $\P_n$ are the $\hat\pi_i(p_n^{-1}(R_{j, n})\cap R_{i, 1})$ for $j=1, \dots, p(n)$.
\end{itemize} 
The fact that all the $y_n$'s satisfy the switching rules implies that the transverse measure is invariant.
\begin{cor}\label{cor}
\begin{itemize}
\item If the dimension of $H_g(S_n, \mr)$ is uniformly bounded by $N$, then  ${\cal M}^t(M)$ contains at most $N$ ergodic measures up to a global rescaling;
\item if furthermore the coefficients of all  the matrices  $A_n$ are uniformly bounded,  then  ${\cal M}^t(M)$ is reduced to a single direction {\em  i.e.} the transverse dynamics is uniquely ergodic. 
\end{itemize}
\end{cor}

\noindent{\it Proof.} The proof is standard and can be found in \cite{GM} in a quite similar situation in the particular case when $g=1$. To prove the first statement we may assume that the dimension of the $H_g(S_n, \mr)$'s is constant and equal to $N$. The set ${\cal M}^t(M)$ is a convex set and its extremal directions coincide with the
set of ergodic measures. Since ${\cal M}^t(M)$is isomorphic to $ \lim_\leftarrow (H_g^+(S_n,\mr),  \tau_{n,\star})$ the convex cone ${\cal M}^t(M)$ is the intersection of the convex nested sets:
$$\M^t(M)\, =\, \cap_{n\geq 0}W_n$$ where
$$W_n\, =\, \tau_{1, \star }\circ \dots\circ \tau_{n-1, \star}(H_g^+(S_n, \mr)).$$ Since each convex cone $W_n$ possesses at most $N$ extremal lines, the limit set $\M^t(M)$ possesses also at most $N$ extremal points and thus at most $N$ ergodic measures.

\noindent In order to prove  the second statement, consider two points $x$ and $y$ in the positive cone of $\mr^N$.
Let $T$ be the largest line segment containing $x$ and $y$ and contained in the positive cone of $\mr^N$. We recall that the hyperbolic distance between $x$ and $y$ is given by:$$Hyp (x, y) \, =\, -\ln{\frac{(m+l)(m+r)}{l.r}},$$ where $m$ is the length of the line segment $[x, y]$ and $l$ and $r$ are the length of the connected components of $T\setminus [x, y]$. Positive matrices contract the hyperbolic distance in the positive cone of $\mr^N$. Since the matrices corresponding to the maps $f_{\star n}$ are uniformly bounded in sizes and entries, this contraction is uniform. Because of this uniform contraction, the set $\M^t(M)$ is one dimensional.

\bigskip
\begin{remark} {\rm  When the Lie group $\mg$ is not amenable, there is no reason a priori that a $\mg$-solenoid $M$ possesses transverse invariant measures. Thus in such a case, when writing $M$ as an inverse limit of $\mg$-branched manifolds:
$$M\, =\, \lim_\leftarrow (S_n, \tau_n),$$
we must have  $H_g^+(S_n,\mr) =0$  for each $n\geq 1$. }\end{remark}

\bigskip
\begin{remark}
{\rm As we have seen, considering the a $\mg$-solenoid as a projective limit of $\mg$-branched manifolds allows us to give a topological characterization of  the invariant measures. It turns out that in the case of $\mr$-solenoid, {\em  i.e.} $\mz$-action free action on the Cantor set, the description in terms of projective limit of graphs, yields a topological characterization of some other dynamical invariants (topological entropy, semi-conjugacy to  rotations etc...) \cite{GM}. It is likely that a similar approach works in   the general case of $\mg$-solenoids. Furthermore, projective limits of graphs can be used to construct $\mr$-solenoid with some nice exotic properties (for example minimal, uniquely ergodic and with infinite topological entropy or with many ergodic invariant measures) and a simple way to construct these examples is to use Toeplitz systems (branches with equal lengths) . At least in the case of $\mr^d$-solenoid, similar constructions can be done by considering $\mr^d$-branched manifolds made of tori with dimension $d$ glued along tori with dimension $d-1$.  } 
\end{remark}

\section{Appendices}
\subsection{ $\me^d$-solenoids}
The Lie group $\me^d$ of  affine orientation preserving isometries of the Euclidean space $\mr^d$ is  the semi-direct product $SO(d) \rtimes \mr^d$. Thus in addition to the $\me^d$-action on a  $\me^d$-solenoid $M$, there exists also a $\mr^d$-(sub) action.  Choose a point $x$ in $M$, we denote by $\omega_{\me^d, \mr^d}(x)$ the closure of the orbit $x+\mr^d$ in $M$.  
\begin{prop} For each $x$ in a  minimal $\me^d$-solenoid $M$,    the $\mr^d$-invariant subset $\omega_{\me^d, \mr^d}(x)$ is  $\mr^d$-minimal {\em  i.e.}all its $\mr^d$-orbits are dense in $\omega_{\me^d, \mr^d}(x)$. 
\end{prop}

\noindent{\it Proof:}   Let $g$ be an element in $\me^d$ and $x$ a point in $M$. Since $(x+\mr^d).g = x.g +\mr^d$, we have: $\omega_{\me^d, \mr^d}(x.g)= \omega_{\me^d, \mr^d}(x).g$.  Any $\mr^d$-action on a compact set possesses a minimal invariant subset, thus there exists $x_0$ in $M$ such that $\omega_{\me^d, \mr^d}(x_0)$ is  $\mr^d$-minimal. It follows that for any $g$ in $\me^d$, $\omega_{\me^d, \mr^d}(x_0.g)$ is either equal to  or completely disjoint from
$\omega_{\me^d, \mr^d}(x_0)$, but in both cases it is also  $\mr^d$-minimal. Since  $M$ is minimal for the $\me^d$-action, any point $x$ in $M$, belongs to a  $\omega_{\me^d, \mr^d}(x_0.g)$ for some $g$ in $\me^d$ and thus 
$\omega_{\me^d, \mr^d}(x)$ is also $\mr^d$-minimal.

Let  $\Gamma_y$ be  the stabilizer of $\omega_{\me^d, \mr^d}(y)$  in $\me^d$  {\em  i.e.} the closed subgroup  of the $\gamma$'s in $SO(d)$ such that $\omega_{\me^d, \mr^d}(y).\gamma = \omega_{\me^d, \mr^d}(y)$. The stabilizers associated with two different points  in $M$ are conjugate.  The proof of the following proposition is immediate:
\begin{prop} \label{iso}Let $M$ be a minimal, $\me^d$-solenoid. Then, for any $y$ in $M$,  the following $2$ assertions are equivalent:
\begin{itemize}
\item [(i)]  $\omega_{\me^d, \mr^d}( y) = \Omega_{\me^d}(y)$; 
\item [(ii)]  $\Gamma_y = SO(d)$. 
\end{itemize}
Furthermore if one of these properties  is satisfied for some $y$ in $M$, then it is also satisfied for any other $y'$ in $M$.
\end{prop}
A  minimal, $\me^d$-solenoid which satisfies one of the properties in Proposition \ref{iso} is called {\it isotropic}. 

Consider now a projective limit of $\me^2$-branched manifolds  associated with a free minimal $\me^2$-solenoid $M = \lim_\leftarrow (S_n, \tau_n)$. The $\me^2$-action on the hull $M$ decomposes itself into two parts when projected on the$\me^2$-branched manifold $S_n$. On one hand, the $SO(2)$-action on   $M$ projects by $p_n$ on the $SO(2)$ action on the $\me^2$-branched manifold $S_n$.  On the other hand, there is of course no $\mr^2$-action on  $S_n$ but for any $y$  in $M$, $\omega_{\me^2, \mr^2}(y)$ projects on a closed subset of $S_{n}$ which is $\Gamma_y$ invariant. The $\me^2$-branched manifold $S_{n}$ is foliated by the projection of the orbits $y+\mr^2$ for $y$ in $M$.  In the particular case when $M$  is  isotropic,  all the projected $\mr^2$-orbits are dense. 
\subsection{ $\mr^d$-solenoids}
In the particular case when the group $\mg =\mr^d$ we have the following result:
\begin{prop} \label{rec} Any  minimal $\mr^d$-solenoid admits a rectangular box decomposition. 
\end{prop}

\noindent{\it Proof:} Consider rectangular boxes that read in  charts of an atlas, $R\times T$, where $R$ is a rectangle with faces parallel to the one of the unit $d$-cube generated by the  canonical basis of $\mr^d$ and $T$ is a totally disconnected set. We call such rectangular boxes  well oriented rectangular boxes The intersection of two well oriented  rectangular boxes is a finite union of well oriented  rectangular boxes and that the closure of the union of two well oriented  rectangular boxes is also the closure of the union of a finite collection of pairwise disjoint well oriented rectangular boxes. It follows that there exists a finite collection rectangular boxes $R_1\times T_1\times, \dots, R_n\times T_n$ pairwise disjoint and whose closures cover the whole solenoid. This yields the following result:
\begin{teo} Any  minimal $\mr^d$-solenoid is orbit equivalent to an $\mr^d$-solenoid  which fibers on the $d$-torus.
\end{teo}

\noindent{\it Proof:}  Consider a minimal $\mr^d$-solenoid $M$. From Proposition \ref{rec}, we can associate with $M$ a rectangular box decomposition $\B$ which by projection along the vertical direction in the boxes gives a $\mr^d$- branched manifold $S_\B$ whose $d$-regions are $d$-dimensional rectangles. Consider the $(1)$-stratum of the stratification of $S_\B$. It is made with a finite number of $1$-dimensional edges parallel to the directions of the vectors of the canonical orthonormal basis and a finite number of vertices.   Assume that all the sizes of these edges are commensurable {\em  i.e.} there exists a positive real $\tau$ such that the $\mr^d$-branched manifold $S$ can be filled with $d$-cubes with size $\tau$ that intersect on their boundary, full faces to full faces. This induces a projection $\pi:S\to \mt^d$ where $\mt^d$ is  the $d$-torus $\mr^d/{\tau\mz^d}$ which commutes with translations. In the general case, the sizes of the rectangles are not commensurable.  However, we can perform a change of the metric $g$ in the leaves that satisfies the following properties:
\begin{itemize}
\item at each point in $S_\B$ the  metric $g$ preserves the orthogonality of basis corresponding in charts to the canonical basis of $\mr^d$;
\item in the chart associated with each rectangular box of the box decomposition,  the metric $g$ is independent on the choice of the point on same vertical;
\item the size of the rectangles are commensurable when computed with the metric $g$. 
\end{itemize}

This new metric defines a homeomorphism between the $\mr^d$-solenoid $S_\B$ and an $\mr^d$-solenoid which possesses a box decomposition in commensurable rectangles and thus fibers on a $d$-torus. As a direct consequence, we recover a result obtained by L. Sadun and R. Williams \cite{SW}:
\begin{cor}Let $X$ be a strongly aperiodic, repetitive, Delone set of $\mr^d$-finite type, then there exists a strongly aperiodic, repetitive, Delone set of  $\mr^d$-finite type $Y$ which is $\mr^d$-orbit equivalent to $X$ and is contained in the lattice $\mz^d$.

\end{cor}

\end{document}